\documentclass{article}

\begin{document}

\title{The ADHM construction of Yang-Mills instantons}
\author{Simon Donaldson}
\date{\today}
\maketitle


\newcommand{\bR}{{\bf R}}
\newcommand{\bC}{{\bf C}}
\newcommand{\bP}{{\bf P}}
\newcommand{\bH}{{\bf H}}
\newcommand{\cO}{{\cal O}}
\newcommand{\ubC}{\underline{{\bf C}}}
\newcommand{\cM}{{\cal M}}
\newcommand{\uE}{\underline{E}}
\newcommand{\uU}{\underline{U}}
\newcommand{\uV}{\underline{V}}
\newcommand{\uW}{\underline{W}}
\newcommand{\uB}{\underline{B}}
\newcommand{\tOmega}{\widetilde{\Omega}}
\newcommand{\cF}{{\cal F}}
\newcommand{\hatT}{\widehat{T}}
\newcommand{\hatvarpi}{\widehat{\varpi}}
\section{The instanton problem}

The mathematical formulation of classical electromagnetism involves  an electromagnetic field $F$ on space-time. This is a skew symmetric tensor $\left(F_{ij}\right)$, or $2$-form. In the classical theory it is useful to write this in terms of an electromagnetic potential $A$:
$$   F_{ij}= \frac{\partial A_{j}}{\partial x_{i}}- \frac{\partial A_{i}}{\partial x_{j}}. $$
In terms of differential forms, $A$ is a $1$-form and $F=dA$. This encodes the fact that $dF=0$, which is a part of the Maxwell equations. The remainder of Maxwell's vacuum equations are, in this notation, $d^{*}F=0$, where $d^{*}$ is the formal adjoint of $d$ defined using the Lorentzian metric on space-time. Explicitly, in terms of a space and time decomposition $x_{0}=t$, the tensor $F$ is written in terms of  electric and magnetic vector fields $E, B$:
$$  F= \sum_{i=1}^{3} E_{i}\  dt dx_{i} + \sum B_{i}\  dx_{j} dx_{k}, $$
where in the second term $(ijk)$ run over cyclic permutations of $(123)$. The equation $dF=0$ is
   \begin{equation}   {\rm curl}\  E = -\frac{\partial B}{\partial t}\ \ \ , \ \ \ {\rm div}\  B=0 \end{equation}
   and the equation $d^{*}F=0$ is
   \begin{equation}   {\rm curl} \ B=  \frac{\partial E}{\partial t}\ \ \ , \ \ \ {\rm div}\  E= 0. \end{equation}

The potential $A$ combines the magnetic vector potential and scalar electric potential.

The potential is not unique, we can change $A$ to $A+d\chi$ for any function $\chi$, so its classical physical meaning is not so clear. In quantum mechanics
it enters as the coupling with electromagnetism. Quantum mechanical formulae such as the Schr\"odinger equation are modified in the presence of an electromagnetic field by replacing the ordinary derivatives $\frac{\partial}{\partial x_{i}}$ acting on wave functions by  
$$ \nabla_{i} = \frac{\partial}{\partial x_{i}} + \sqrt{-1}\  A_{i}. $$
  This is  all clarified by the differential  geometric language of vector bundles and connections. The wave function is a section of a complex Hermitian line bundle over space-time,  the potential is a connection on this line bundle and  (3) is the usual formula for the covariant derivative of a section.
The ambiguity in the connection now appears as the choice of trivialisation of the line bundle. 
More generally in  Yang-Mills theory, beginning in 1954, this discussion is extended to consider  vector bundles with connections of higher rank, or more generally still bundles with structure group some Lie group $G$. Then the \lq\lq potential'' is the connection $1$-form,  with  values in   the Lie algebra of $G$. We will restrict attention to matrix groups $G$ so one can think of a connection as a covariant derivative $\nabla_{i}=\partial_{i}+ A_{i}$ but now acting on vector-valued functions. The generalisation of the electro-magnetic field  is the curvature
$$     F_{ij}= \frac{\partial A_{j}}{\partial x_{i}}- \frac{\partial A_{i}}{\partial x_{j}} + [A_{i}, A_{j}],  $$
which appears as  the commutator $[\nabla_{i},\nabla_{j}]$. 

(See Yau's lecture \cite{kn:Yau} in this series for more on the  evolution of these ideas.) The analogues of  Maxwell's equations are the Yang-Mills equations. These are obtained from a Lagrangian
$$ {\cal E}(A)=  \int \vert F\vert^{2}. $$
Here the integrand $\vert F\vert^{2}$ is defined using the metric tensor on space time and an inner product on the Lie algebra of the structure group $G$. The Euler-Lagrange equation $\delta {\cal E}=0$ is the Yang-Mills equation
 $$ d^{*}_{A}F_{A}=0, $$
 where $d^{*}_{A}$ is the coupled formal adjoint. This is a nonlinear second order PDE for the connection $A$.

So far we have been working on Lorentzian space time, but now we consider the case of Euclidean $\bR^{4}$, which is relevant in Quantum Field Theory. The formalism is the same but the Yang-Mills equations become elliptic (when interpreted modulo the gauge freedom) rather than hyperbolic. All of this can be done in any dimension but there are special features in dimension $4$, one being  that the equations are conformally invariant. The solutions of the Yang-Mills equation on $\bR^{4}$ which extend to the conformal
 compactification $S^{4}$ are exactly those with finite energy ${\cal E}$.  Another special feature in dimension $4$ is that  there are first order \lq\lq instanton'' equations which imply the second order Yang-Mills equations. This is analogous to the relation in two dimensions between the Cauchy-Riemann and Laplace equations for complex-valued functions.  The curvature $F$ of a connection is a bundle valued $2$-form and the instanton equation is $F= - *F$ where $*:\Lambda^{2}\rightarrow \Lambda^{2} $ is the Hodge $*$-operator, which is also conformally invariant. These instantons, when they exist, are absolute minima of the energy on a given $U(r)$ bundle $E\rightarrow S^{4}$. In fact the energy is the topological invariant $8\pi^{2} c_{2}(E)$, where $c_{2}(E)$ denotes the second Chern class, evaluated on the fundamental class of the sphere. Thus an instanton can only exist if $c_{2}\geq 0$, for negative $c_{2}$ we get an equivalent discussion by changing the orientation. If we go back to $\bR^{4}$ and choose a  \lq\lq space-time'' decomposition $\bR\times \bR^{3}$ (but now with a positive metric) the curvature can be written as a pair of \lq\lq electric'' and \lq\lq magnetic'' components $(E,B)$ and the instanton equation is $E=B$. (In Lorentzian space the equations (1),(2) are symmetric between $E,B$ up to sign changes, but in Euclidean signature there is an exact symmetry.)

We have now set the scene for the problem: {\em for each $r\geq 2$ and $k\geq 1$ describe all the solutions of the instanton equation over $S^{4}$ on a $U(r)$ bundle with $c_{2}=k$.}

\section{The ADHM construction}

This problem was solved by 
 Atiyah, Drinfeld, Hitchin and Manin (ADHM) in 1977.
 In his commentary on Volume 5 of his collected works \cite{kn:Atiyah}, Atiyah wrote:

\

{\em \dots    with the help of Nigel Hitchin, I finally saw how Horrocks' method gave a very satisfactory and explicit solution to the problem. I remember our final discussion one morning when we had just seen how to fit together the last pieces of the puzzle. We broke off for lunch feeling very pleased with ourselves. On our return, I found a letter from Manin (whom I had earlier corresponded with on this subject) outlining essentially the same solution to the problem and saying \lq\lq no doubt you have already realised this''! We replied at once and proposed that we should submit a short note from the four of us.}

\

This short note appeared as \cite{kn:ADHM}. Atiyah gave a detailed exposition of the proof and the background in his Pisa lectures \cite{kn:AtiyahP}.

The ADHM construction can be explained in relatively elementary terms. For fixed $k$ and $r$ consider  a family of linear maps $\lambda_{x}:\bC^{2k+r}\rightarrow \bC^{2k}$ parametrised by $x\in \bR^{4}$ of the form 
$$  \lambda_{x}= \sum_{i=0}^{4} x_{i} L_{i}  + M, $$
where $L_{0},L_{2}, L_{3}, M$ are $2k\times (2k+r)$ matrices. Suppose that
$\lambda_{x}$ is surjective for all $x\in \bR^{4}$. Then $\uE={\rm ker} \lambda$ is a rank $r$ bundle over $\bR^{4}$. Suppose also that the image of the $L_{i}$ span $\bC^{2k}$. Then  the bundle has a natural extension to the $4$-sphere with $c_{2}=k$. The statement is that if the matrix data $L_{i},M$ satisfies a system of quadratic equations (given explicitly in (8),(9) below) then the induced connection on $\uE$ is an instanton connection and all arise in this way. Moreover the matrix data associated to an instanton is unique up to the obvious action of $U(2k)\times U(2k+r)$. In short, the solution of the instanton equations (which are nonlinear first order PDE) is reduced to solving algebraic equations for the matrix data.

The induced connection referred to above is the same construction as in classical differential geometry for the tangent bundle of a submanifold of Euclidean space. Writing $\ubC^{2k+r}$ for the trivial bundle with fibre $\bC^{2k+r}$ we have an inclusion $\iota : \uE\rightarrow \ubC^{2k+r}$ and an orthogonal projection map $\pi: \ubC^{2k+r}\rightarrow \uE$, using the standard Hermitian metric on $\bC^{2k+r}$. Then the covariant derivative of the induced connection on $\uE$ is defined by
\begin{equation}    \nabla s = \pi\circ \nabla_{{\rm flat}} \iota, \end{equation}
where $\nabla_{{\rm flat}}$ is the usual derivative on sections of the trivial bundle $\ubC^{2k+r}$. 

 While the construction  can be described in elementary terms, as above, the  fact that it gives the general solution is much deeper. The ADHM work   is significant as one of the first applications of sophisticated modern geometry to physics, playing a large part in opening up a dialogue which has of course flourished mightily in the near half century since. It is also a beautiful piece of mathematics which can be approached from many directions.

\section{Twistor space and the Ward correspondence}

The approach of Atiyah, Drinfeld, Hitchin and Manin to this problem went through the Penrose twistor theory and the Ward correspondence between instantons and holomorphic vector bundles, which we will now review.

Let $Z$ be a complex manifold and $W\subset Z$  a compact complex submanifold
with normal bundle $N$. The theory of Kodaira  describes the small deformations of  $W$ in terms of the cohomology groups $H^{0}(N)$ and $ H^{1}(N)$. If $H^{1}(N)=0$ then the moduli space of deformations is a complex manifold  ${\cal M}$,with tangent space at $W$ equal to the space of holomorphic sections $H^{0}(N)$. Now suppose that $Z$ has complex dimension $3$ and that the submanifold is a \lq\lq line'' $L$---
an embedding of the Riemann sphere $\bC\bP^{1}$. Vector bundles over $\bC\bP^{1}$ are all equivalent to direct sums of line bundles $\cO(p)$, the tensor powers of the Hopf line bundle $\cO(1)$. Suppose that the normal bundle of $L$ is isomorphic to $\cO(1)\oplus \cO(1)$. Then $H^{1}(N)=0$ and $H^{0}(N)$ is $4$-dimensional so we obtain a $4$-dimensional complex manifold ${\cal M}$ of lines in $Z$ whose points represent deformations of $L$ with normal bundle $\cO(1)\oplus \cO(1)$.

In this situation there is an induced holomorphic conformal structure on ${\cal M}$.  For if we write $N= S_{-}\otimes \cO(1)$ for a two dimensional complex vector space $S_{-}$ then $H^{0}(N)= S_{-}\otimes S_{+}$ where
$S_{+}$ is the $2$-dimensional complex vector space of sections of $\cO(1)$ over $L$. The tensor product of the skew-symmetric maps $S_{\pm}\otimes S_{\pm}\rightarrow \Lambda^{2} S_{\pm}$ is a symmetric map
$$   T\cM\otimes T\cM\rightarrow \Lambda^{2}S_{+}\otimes \Lambda^{2} S_{-}$$
which gives the conformal structure (that is, it becomes a nonsingular quadratic form when one fixes an identification of the line $\Lambda^{2}S_{+}\otimes \Lambda^{2} S_{-}$ with $\bC$).

The geometric meaning of this construction is that the null cone in the tangent space of $\cM$ at a line $L$ is the set  of infinitesimal deformations of $L$ which intersect $L$. For each point $p\in Z$ we get a submanifold 
$\Sigma_{p}\subset \cM$ consisting of lines through $p$ and this is isotropic, {\it i.e.} the conformal structure vanishes on the tangent spaces of $\Sigma_{p}$.

We can now consider the \lq\lq instanton equations''  in a holomorphic setting, for holomorphic connections on a bundle over $\cM$. Let $Q$ be a nondegenerate quadratic form on $\bC^{4}$. The complex linear $*$-operator on $2$-forms  is defined  just as in the real case and hence the decomposition of $\Lambda^{2}\bC^{4}$ into self-dual and anti-self-dual components. A $2$-form is $\pm$ self-dual if and only if its restriction to each isotropic subspace in $\bC^{4}$ is zero. Applying this to the
 tangent spaces of $\cM$ we see that a connection is an instanton if and only its restriction to each of the submanifolds $\Sigma_{p}$ is flat.

Suppose that $E$ is a holomorphic vector bundle over $Z$ which is holomorphically trivial on each line. Then we define a bundle $\uE$ on $\cM$ with fibre over
a line $L$ equal to the holomorphic sections of $E$ over $L$. The basic point is that $\uE$ has a natural connection. To see this one can use the theory of formal neighbourhoods. The first formal neighbourhood $L_{(1)}$ of a line $L$ is the sheaf $\cO_{Z}/ {\cal I}^{2}$ where ${\cal I}\subset \cO_{Z}$ is the ideal sheaf of functions on $Z$ vanishing on $L$. The obstruction to extending a trivialisation of $E$ over $L$ to $L_{(1)}$ lies in $H^{1}(L; N^{*}\otimes E\vert_{L})$ and if this vanishes any two extensions differ by an element of $H^{0}(L; N^{*}\otimes E\vert_{L})$. Since the cohomology groups
 $H^{p}(\cO(-1))$ vanish there is a unique extension of a trivialisation and this defines a connection on $\uE$. For each $p$ in $Z$ we get a canonical trivialisation of $\uE$ over $\Sigma_{p}$ by identifying sections with their values at $p$. These trivialisations are compatible with the connection, so the connection is flat over the $\Sigma_{p}$ and hence is an instanton. Conversely, starting with an instanton $(\uE,\nabla)$ over $\cM$ we get a holomorphic bundle $E$ over $Z$ with fibre over $p$ equal to the covariant constant sections of $\uE$ over $\Sigma_{p}$. These constructions are mutually inverse and give the Ward correspondence between instantons on $\cM$ and holomorphic bundles over $Z$ which are trivial on all lines.

To get back to real $4$-manifolds one considers a complex manifold $Z$ with a \lq\lq real'' structure: an antiholomorphic involution $\sigma:Z\rightarrow Z$. Then the space of \lq\lq real lines''---{\it i.e.} the lines preserved by $\sigma$ is a real form $M$ of $\cM$. There are two cases: either $\sigma$ acts on the real lines as complex conjugation, with fixed point set a circle, or as an antipodal map, with no fixed points. The first case corresponds to a conformal structure on $M$ of signature $(2,2)$ and the second case to a Riemannian conformal structure and this latter case is the one relevant to the problem  introduced in the previous section.

With this general background in place we focus on the twistor space $Z=\bC\bP^{3}$. The space of lines is the complex Grassmannian ${\rm Gr}(2,4)$ which can be viewed as  a quadric in $\bC\bP^{5}$. If we identify $\bC^{4}$ with the quaternionic vector space $\bH^{2}$ we get a fixed-point-free real structure
$\sigma$ induced by multiplication by the quaternion $J$. The real lines are the fibres of the natural map
  $$\pi: \bC\bP^{3}\rightarrow \bH\bP^{1}= S^{4}, $$
and instantons on $S^{4}$ correspond to holomorphic bundles on $\bC\bP^{3}$ trivial on all real lines. An instanton on $S^{4}$ with structure group $U(r)$ corresponds to a holomorphic bundle $E$, trivial on all real lines, together with an isomorphism
$\overline{E}^{*}\cong \sigma^{*}(E)$. 

Another approach to the theory in the Riemannian case, following Atiyah, Hitchin, Singer \cite{kn:AHS} is to obtain the twistor space of an \lq\lq anti-self-dual'' Riemannian 4-manifold $M$ as the projectivization of the positive spin bundle. A point in $Z$ is a pair consisting of a point  in $M$ and a compatible complex structure on the tangent space at that point. In this approach the bundle $E$ is just the pull-back $\pi^{*}(\uE)$ and the Newlander-Nirenberg theorem is used to produce a holomorphic structure. In the holomorphic approach discussed above one would first show, by elliptic PDE theory, that an instanton is real analytic in a suitable gauge and thus extends to a complexification of the $4$-manifold.

\section{Construction of bundles over $\bC\bP^{3}$.}

For a bundle $E$ over $\bC\bP^{n}$ we  use the standard notation $E\otimes \cO(p)= E(p)$, where $\cO(p)$ is the $p$th. power of the Hopf line bundle $\cO(1)$.

Let $U,V,W$ be complex vector spaces and suppose that we have bundle maps
\begin{equation}   \uU(-1) \stackrel{a}{\rightarrow} \uV \stackrel{b}{\rightarrow} \uW(1), \end{equation}
over $\bC\bP^{3}$
with $a$ injective, $b$ surjective and $b\circ a=0$. (Here $\uU$ denotes the trivial bundle with fibre $U$.) Then we get a holomorphic bundle $E$ as the cohomology ${\rm Ker}\ b/{\rm Im}\ a$. The data (4) is called a \lq\lq monad'' and the construction was introduced by Horrocks in 1964 \cite{kn:Horrocks}. 
Explicitly, with homogeneous co-ordinates $Z_{i}$ and with chosen bases for the three vector spaces we write
$a= \sum A_{i} Z_{i}\ ,\ b=\sum B_{i} Z_{i}$ for matrices $A_{i}, B_{i}$. A monad is given  by a solution of the matrix equations
\begin{equation}     B_{j}A_{i}+ B_{i} A_{j}= 0 \ \ \ \ (i,j=1,\dots,4) \end{equation}

Following Horrocks, the monad construction was studied by Barth  and Barth and Hulek \cite{kn:BarthHulek} who proved:

\

{\bf Theorem}

\

{\em Let $E$ be a rank $r$ holomorphic bundle over $\bC\bP^{3}$ with $c_{2}=k$. Suppose that $E$ is trivial on some line and satisfies the vanishing conditions
$$  \ H^{0}(E)=0\ \ ,\ \  H^{0}(E^{*})=0, \ \ H^{1}(E(-2))=0 \ \ , \ \ H^{1}(E^{*}(-2))=0. $$
Then $E$ arises from a monad with ${\rm dim} \ U={\rm dim}\ W=k$ and ${\rm dim} V= 2k+r$. This monad is unique up to the action of $GL(U)\times GL(V)\times GL(W)$.}

\

(Note: In fact Barth and Hulek considered symplectic or orthogonal bundles $E$, with $E\cong E^{*}$.)

\

Let $L$ be a line in $\bC\bP^{3}$ and $p,q$ be distinct points on $L$. Given a monad we have linear maps $a_{p}, a_{q}: U\rightarrow V$ and $ b_{p}, b_{q}: V\rightarrow W$. One finds that the cohomology bundle $E$ is trivial on $L$ if and only if $b_{q}\circ a_{p}$ is an isomorphism from $U$ to $W$, which is equivalent to the same for $b_{p}\circ a_{q}$ by the equations (6). In this case the fibres of the restriction of $E$ to $L$ can be written both as a subspace
   $$  {\rm Ker}\ b_{p} \cap {\rm Ker}\  b_{q}$$ of $V$ and as a quotient
   $$  V/({\rm Im}\ a_{p} + {\rm Im}\  a_{q}).$$
   Thus the associated bundle $\uE$ over $S^{4}$ arises as both a subbundle and a quotient of the trivial bundle with fibre $V$. One finds that the instanton connection is that given by the standard construction (3). 

Given the Barth-Hulek theorem above we  get a complete description of instantons on $S^{4}$ if it can be shown that the corresponding bundles over $\bC\bP^{3}$ satisfy the vanishing conditions.

The condition that $H^{0}(E)=0$ is relatively easy to understand. From the construction, a non-trivial holomorphic section of $E$ goes over to a covariant constant section of the bundle $\uE$ over $S^{4}$ which means that the structure group reduces to $U(r-1)$. In view of the symmetry between $E, E^{*}$ in the hypotheses of the theorem the essential condition is $H^{1}(E(-2))=0$. The reason that this holds goes back to the origins of Penrose's twistor theory. 

Let $\Omega \subset \bR^{4}$ be an open set and $\tOmega\subset \bC\bP^{3}$ be the union of the corresponding lines. Let $\chi$ be a sheaf cohomology class in $H^{1}(\tOmega, \cO(-2))$. For a line $L$ we have $H^{1}(L;\cO(-2))= \bC$ so the restriction of $\chi$ to lines gives a function $F_{\chi}$ on $\Omega$. More precisely, there is  a natural isomorphism of $H^{1}(L, K_{L})$ with $\bC$ where $K_{L}$ is the canonical line bundle of $L$. So to define $F_{\chi}$ we need to specify an isomorphism of $K_{L}$ with $\cO(-2)\vert_{L}$.
If $L=\bP(S)$ this is the same as a trivialisation of $\Lambda^{2}S$. Let $L_{\infty}$ be the line corresponding to the point at infinity in $ S^{4}$ and write $L_{\infty}= \bP(S_{\infty})$. Then for a line $L=\bP(S)$ corresponding to a point of $\bR^{4}$ we have $\bC^{4}= S\oplus S_{\infty}$ which gives the trivialisation of $\Lambda^{2}S$.

Penrose showed that the function $F_{\chi}$ satisfies the Laplace equation on $U\subset \bR^{4}$ and moreover that this gives an equivalence between $H^{1}(\tOmega, \cO(-2))$ and the solutions of the Laplace equation. In terms of a C\u{e}ch representation of $\chi$ by a holomorphic function of three complex variables the function $F_{\chi}$ is given by a contour integral formula of a classical nature---see Hitchin's article \cite{kn:Hitchin2} in this series.
Penrose's theory interprets other sheaf cohomology groups on twistor space in terms of other linear field equations. In particular this  correspondence for $H^{1}(\cO(-2)$ can be related to the instanton construction we discussed above. Suppose that $\Omega$ is simply connected and consider the instanton equation on $\Omega$ for the structure group $\bC^{*}$, {\it i.e.} for line bundles. The Ward correspondence tells us that these are determined by line bundles over $\tOmega$ which in this case correspond to classes in $H^{1}(\tOmega, \cO)$. Take a pair of planes in $\bC\bP^{3}$ containing the line at infinity and let $\Sigma_{+},\Sigma_{-}$ be their intersections with $\tOmega$. Then classes in $H^{1}(\tOmega, \cO(-2))$ correspond to holomorphic line bundles on $\tOmega$ trivialised over $\Sigma_{+}$ and $\Sigma_{-}$. A complexification of $\Omega$ in the space of lines can be identified with $\Sigma_{+}\times \Sigma_{-}$. We get complex co-ordinates $(z_{1},z_{2}, w_{1}, w_{2})$ in which the metric is $dz_{1}dw_{1} + dz_{2}dw_{2}$. The instanton equations for a connection $1$-form
$$  \alpha_{1} dz_{1} + \alpha_{2} dz_{2}+ \beta_{1} dw_{1}+ \beta_{2}dw_{2}, $$ are
\begin{equation}  \frac{\partial \alpha_{1}}{\partial z_{2}} =  \frac{\partial \alpha_{2}}{\partial z_{1}}\ \ ,\ \    \frac{\partial \beta_{1}}{\partial w_{2}} =  \frac{\partial \beta_{2}}{\partial w_{1}}\ \ ,\ \  \frac{\partial \alpha_{1}}{\partial w_{1}} + \frac{\partial \alpha_{2}}{\partial w_{2}} = \frac{\partial \beta_{1}}{\partial z_{1}} +\frac{\partial \beta_{2}}{\partial z_{2}} . \end{equation}

 Following though the Ward correspondence, one find that a holomorphic line bundle on $\tOmega$ trivialised over $\Sigma_{\pm}$ gives a connection $1$-form determined by two functions $f,g$ in the shape
$$  \alpha_{i}= \frac{\partial f}{\partial z_{i}}, \beta_{i}= \frac{\partial g}{\partial w_{i}}. $$

The last equation in (6) then shows that the function $F=f-g$ satisfies
$$   \frac{\partial^{2} F}{\partial z_{1}\partial w_{1}}+ \frac{\partial^{2} F}{\partial z_{2}\partial w_{2}}=0 $$
which is the Laplace equation for the metric $dz_{1}dw_{1} + dz_{2}dw_{2}$.

The result from the Penrose theory that we require extends that discussed above in two ways. First, we need to consider the whole of $S^{4}$ rather than an open subset in $\bR^{4}$ and, second, the cohomology group $H^{1}(E(-2))$ in question involves the bundle $E$. The conclusion is that the group corresponds to sections $s$  of the instanton bundle $\uE$ over $S^{4}$ satisfying the equation
$$\nabla^{*}\nabla s + \frac{R}{8} s =0, $$
where $R$ is the scalar curvature of the round metric on $S^{4}$. This is the coupled, conformally-invariant, Laplace equation. Since the scalar curvature is positive, integration by parts shows that the only solution is $s=0$, which translates back to the statement $H^{1}(E(-2))=0$. 

\section{The Beilinson spectral sequence}

The quickest way to establish the Barth-Hulek Theorem  uses a construction of Beilinson \cite{kn:Beilinson} which appeared a little after the ADHM work.

Let $z_{1}, z_{2}, z_{3}$ be standard co-ordinates on $\bC^{3}$ and consider the vector field 
$$  v= \sum z_{i} \frac{\partial}{\partial z_{i}}, $$
which has a single zero at the origin. This has a holomorphic extension to a vector field on $\bC\bP^{3}$ vanishing on the plane at infinity. Thus it defines a section, which we also call $v$, of $T(-1)$, where $T$ is the tangent bundle of $\bC\bP^{3}$. This section has  a single zero at the origin $0\in \bC^{3}\subset \bC\bP^{3}$.

Contraction with a non-zero vector $e$ in a vector space $K$ defines an exact sequence
$$   \dots \Lambda^{r}K^{*}\rightarrow \Lambda^{r-1}K^{*}\dots \rightarrow K^{*}\rightarrow \bC. $$
In this way, the section $v$ gives a {\it Koszul complex}
\begin{equation}  \Omega^{3}(3)\rightarrow\Omega^{2}(2) \rightarrow \Omega^{1}(1) \rightarrow \cO \rightarrow \cO_{0} \end{equation}
which is an exact sequence of sheaves. Here the $\Omega^{r}$ denote the exterior powers of the cotangent bundle. The top power $\Omega^{3}$ is the canonical bundle $\cO(-4)$ so the first term is $\cO(-1)$. The final term $\cO_{0}$ is the skyscraper sheaf at $0$. 
Let $V$ be any vector bundle over $\bC\bP^{3}$ and take the tenor product of (8) with $V$ to get an exact sequence of sheaves, say,
$$   \cF(3)\rightarrow \cF(2)\rightarrow \cF(1)\rightarrow \cF_{0}\rightarrow V_{0}. $$
In this general situation there is a \lq\lq hypercohomology spectral sequence'' with $E_{1}$ page
  $$ E^{p,q}_{1}= H^{q}(\cF(-p)) $$
  abutting to the vector space $V_{0}$,  the fibre over $V$ at $0$, in degree $0$. (See 3.5 on \cite{kn:GH}, for example.)

In our situation, take $V=E(-2)$ where $E$ a bundle satisfying the hypotheses of the Barth-Hulek theorem. The $E_{1}$ page has rows
$$  H^{q}(E(-3))\rightarrow H^{q}(E\otimes \Omega^{2})\rightarrow H^{q}(E\otimes \Omega^{1}(-1))\rightarrow H^{q}(E(-2)), $$for $0\leq q\leq 3$.

We have the exact Euler sequence
$$  0\rightarrow \cO(-1)\rightarrow \ubC^{4}\rightarrow \Omega^{2}(3)\rightarrow 0 $$
and its dual
$$  0\rightarrow \Omega^{1}(1)\rightarrow \ubC^{4}\rightarrow \cO(1)\rightarrow 0. $$
(In the first sequence we have used the fact that $T\cong \Omega^{2}(4)$).
Using these, Serre duality and the cohomology vanishing assumptions, we find that all the
$H^{0}$ and $H^{3}$ terms vanish. In the $H^{1}$ row we know that $H^{1}(E(-2))=0$ and the dual Euler sequence, tensored with $E(-2)$ shows that $H^{1}(E\otimes \Omega^{1}(-1))=0$. Let $P$ be a plane in $\bC\bP^{3}$ containing a line on which $E$ is trivial. We have a  long exact cohomology sequence which contains
$$  H^{0}(E(-2))\vert_{P}\rightarrow H^{1}(E(-3)\rightarrow H^{1}(E(-2)). $$

The condition of being trivial on a line is open, so $E$ is trivial on the generic line in $P$ and it follows that $H^{0}(E(-2)\vert_{P})=0$. Then the exact sequence shows that $H^{1}(E(-3))=0$. Since the  spectral sequence abuts to a non-zero term only in total degree $0$ we must have $H^{1}(E\otimes \Omega^{2})=0$. So, in sum, all the $H^{1}$ terms vanish. Finally, $H^{2}(E(-2)$ is the Serre dual of $H^{1}(E^{*}(-2)$ and so this vanishes. Thus the $E_{1}$ page can be written as
$$  H^{2}(E(-3))\stackrel{a}{\rightarrow} H^{2}(E\otimes \Omega^{2})
\stackrel{b}{\rightarrow} H^{2}(E\otimes \Omega^{1}(-1)), $$ and it follows from the general hypercohomology theory that $a$ is injective, $b$ is surjective and that there is a natural isomorphism of the cohomology ${\rm Ker}\ b/{\rm Im}\ a$ with the fibre of $E(-2)$ at $0$. 

Of course there is nothing special about the point $0\in \bC\bP^{3}$ and for any point $x\in \bC\bP^{3}$ we get a similar description of the fibre of $E(-2)$ over $x$. When we keep track of the $x$ dependence we find that this is precisely a description of $E$ as a monad (5), with $U=H^{2}(E(-3)), V=H^{2}(E\otimes \Omega^{2}), W=H^{2}(E\otimes \Omega^{1}(-1))$. This  completes the proof of the Theorem. (Another way of expressing the last step is to use a resolution of the diagonal in $\bC\bP^{3}\times \bC\bP^{3}$.)

There is a subtlety here involving duality. If we have a monad description (4) of $E$ we get a description of $E^{*}$ with the monad
$$    W^{*}(-1)\stackrel{b^{T}}{\rightarrow} V^{*}\stackrel{a^{T}}{\rightarrow}
 U^{*}(1) . $$
 But it is not immediately obvious that the descriptions above for $E,E^{*}$ are related in this way. To explain this, recall that the hyperplane class in $H^{2}(\bC\bP^{3}, \bC)$ can be regarded as an element of $H^{1}(\Omega^{1})$. The cup products with this class define maps
$$  \chi_{1}: H^{1}(E\otimes \Omega^{1})\rightarrow H^{2}(E\otimes \Omega^{2}), $$
$$   \chi_{2}: H^{1}(E(-1))\rightarrow H^{2}(E\otimes \Omega^{1}(-1)), $$
and it is an exercise to show that these are isomorphims. Under these isomorphisms
the map $b$ corresponds to a map
$$  H^{1}(E\otimes \Omega^{1})\rightarrow H^{1}(E(-1)). $$
These vector spaces are the Serre duals of $H^{2}(E^{*}\otimes \Omega^{2}), H^{2}(E^{*}(-3))$ and the map is the transpose of the map $a$ for $E^{*}$. 
\

The only twistor spaces which are complex projective manifolds are $\bC\bP^{3}$ and the three dimensional flag manifold. The latter corresponds to the $4$-manifold $\bC\bP^{2}$ with its Fubini-Study metric and reverse of the standard orientation.
A description of instantons on $\bC\bP^{2}$ following the ADHM pattern was found by Buchdahl \cite{kn:Buchdahl}.

\section{Explicit matrix description and the Euclidean approach}

Putting everything together, we find that the general solution of the instanton equation can be obtained from the construction described in Section 2 with a family of linear maps $\lambda_{x}= \sum L_{i} x_{i} + M: \bC^{2k+r}\rightarrow \bC^{2k}$.  We write
$\bC^{2k+r}= \bC^{k}\oplus \bC^{k}\oplus \bC^{r}$ and use block notation. Then we can normalise so that:

$$  L_{0}= \left(\begin{array}{ccc} 1&0&0\\0&1&0\end{array}\right); $$

$$  L_{1}= \left(\begin{array}{ccc} -i&0&0\\0&i&0\end{array}\right); $$

$$  L_{2}= \left(\begin{array}{ccc} 0&1&0\\-1&0&0\end{array}\right); $$
$$  L_{3}= \left(\begin{array}{ccc} 0&-i&0\\-i&0&0\end{array}\right); $$
and
$$  M= \left(\begin{array}{ccc} \alpha^{*}_{1}&\alpha^{*}_{2}&P^{*}\\-\alpha_{2}&\alpha_{1} &Q\end{array}\right). $$
This form incorporates the reality conditions. We have a pair of $k\times k$  matrices $\alpha_{1}$ and $ \alpha_{2}$,  a $k\times r$  matrix $P$ and an $r\times k$ matrix $Q$, all with complex co-efficients. The algebraic equations which have to be solved are 
\begin{equation}   [\alpha_{1}, \alpha_{2}] +PQ=0 \end{equation}
\begin{equation} [\alpha_{1},\alpha_{1}^{*}]+[\alpha_{2}, \alpha_{2}^{*}] - PP^{*}+Q^{*} Q =0. \end{equation}
We also have the   nondegeneracy condition, that for each $x\in \bR^{4}$ the map $\lambda_{x}$ is surjective. The moduli space $M_{k,r}$ of instantons is the quotient of the set of matrices $(\alpha_{1}, \alpha_{2}, P, Q)$ satisfying these conditions by the action of $U(k)\times U(r)$.

There are many different notations which can be used to write these equations. The one above is adapted to a choice of complex structure on $\bR^{4}$ {\it i.e.} we write $\bR^{4}=\bC^{2}$ and $z_{1}= x_{0} + \sqrt{-1} x_{1}, z_{2}= x_{2}+\sqrt{-1} x_{3}$ but one can also use spinorial or  quaternionic notations.

\

A few years after the ADHM work, Corrigan and Goddard showed how to derive the construction by direct calculation, without going through twistor theory \cite{kn:CG}. The essential object is the space $H$ of $L^{2}$ solutions over $\bR^{4}$ of the Dirac equation for negative spinors coupled to the instanton connection. This has dimension $k$ by the index theorem.  For $\psi, \psi '\in H$ and each co-ordinate
$x_{i}$ define
$$   T_{i}(\psi,\psi') = \sqrt{-1} \int_{\bR^{4}} x_{i}\  (\psi,\psi')  $$
Standard theory shows that $\psi_{1}, \psi_{2}$ decay as $O(r^{-3})$, so the integrals  converge. Using the $L^{2}$ metric on $H$ the $T_{i}$ can be regarded as skew-adjoint operators on $H$. Then the matrix data $\alpha_{1}, \alpha_{2}$ is given by
$$  \alpha_{1}= T_{0} + \sqrt{-1} T_{1} \ \ \ ,\ \ \ \alpha_{2}=T_{2}+ \sqrt{-1}T_{3}. $$
The remaining data $P,Q$ is obtained from the asymptotics of elements of $H$. With suitable conformal weights, the  Dirac operator is conformally invariant  so $H$ can also be viewed as solutions of the coupled Dirac equation over $S^{4}$. Thus we have an evaluation map
$$  {\rm ev}: H\rightarrow \uE_{\infty}\otimes S^{-}_{\infty}$$
where $\uE_{\infty}$ is the fibre at infinity  of the bundle $\uE$ carrying  the instanton connection and $S_{-}$ is the negative spin space at $\infty$.
From the Euclidean point of view, ${\rm ev}(\psi)$ defines the $O(r^{-3})$ term of $\psi$ as $r\rightarrow \infty$.  Then with respect to a basis of the two dimensional space $S^{-}_{\infty}$
$$ {\rm ev}= P\oplus Q^{*}: H\rightarrow E_{\infty}\oplus E_{\infty}. $$
Then our family of linear maps $\lambda_{x}$ is  a family of natural maps

$$\lambda_{x}: H\oplus H\oplus E_{\infty} \rightarrow H\oplus H , $$

which can be defined using the geometry and analysis of the Dirac equation on $\bR^{4}$. The connection with the complex geometry \lq\lq monad'' approach and spectral sequences is discussed in \cite{kn:DK}. As one step in that direction, the Serre dual of the space $H^{2}(E(-3))$ appearing in the monad is $H^{1}(E^{*}(-1))$ and this corresponds to the coupled harmonic spinors  via the Penrose theory. 

\section{Mukai and Nahm tranforms}

Around the same time as the ADHM work, but coming from a completely different direction, Mukai introduced a construction which, it was later realised, has many formal similarities. Mukai considered a complex torus $T$ and its dual $\widehat{T}$. Thus a point of $\xi$ of $\widehat{T}$ correspond to a flat holomorphic line bundle $L_{\xi}$ over $T$. Let $E$ be a holomorphic vector bundle over $T$. For each $p\leq {\rm dim} T$ we have a family of vector spaces $H^{p}(T;E\otimes L_{\xi})$ parametrised by $\widehat{T}$. In general the dimension of these vector spaces can jump but in a case when for all $\xi$ there is cohomology in a single dimension $p$ this family defines a vector bundle $\hat{E}$ over $\hatT$. (In general,  Mukai's construction operates in derived categories of coherent sheaves.) The relevant case for us is when ${\rm dim}_{\bC}\ T=2$ and we have a bundle $E$ such that the $H^{p}(E\otimes L_{\xi})$ are zero for $p=0,2$. 

 The striking fact is that this \lq\lq Mukai transform''
is symmetric between $T,\hat{T}$. By the standard theory of the Poincar\'e bundle, the dual of $\hatT$ is $T$. If we start with a bundle $E\rightarrow T$ as above and assume that the bundle $\hat{E}\rightarrow \hatT$ satisfies the same vanishing conditions then the double transform recovers the bundle $E$, up to pulling back by the map $x\mapsto -x$ from $T$ to $T$. There is a close analogy with the Fourier transform. 

One way to relate the Mukai construction to instantons goes via twistor theory.
Let $M=\bR^{4}/\Lambda$ be a  flat {\it Riemannian} 4-torus. The  complex structures on  $\bR^{4}$ compatible with the given Euclidean structure are parametrised by a $2$-sphere $S^{2}$ and the twistor space $Z$ is just $S^{2}\times M$, as a smooth manifold. The complex structure on $Z$ is such that the projection
$\Pi:Z\rightarrow  S^{2}$ is holomorphic and the fibres are complex tori corresponding to the different compatible complex structures on $M$. Let $\hat{M}$ be the dual Riemannian torus with its own twistor space $\hat{Z}$ and  fibration $\hat{\Pi}:\hat{Z}\rightarrow S^{2}$. The fibres of $\Pi$ and $\hat{\Pi}$ over a  point of $S^{2}$ are dual complex tori. Then we can make a fibrewise Mukai transform, which takes suitable holomorphic bundles on $Z$ to holomorphic bundles on $\hat{Z}$. In turn, via the Ward correspondence, this gives a transform from instantons over $M$ to instantons over $\hat{M}$. Expressed directly in Riemannian terms, starting with an instanton connection $\nabla$ on a rank r complex vector bundle $E\rightarrow M$ we define a family of vector spaces parametrised by $\hat{M}$ as the kernels of the coupled Dirac operator on sections of $S^{-}\otimes E\otimes L_{\xi}$. If the connection $\nabla$ is irreducible these form a vector bundle $\hat{E}$ over $\hat{M}$ which has a connection $\hat{\nabla}$ induced by the a version of the  projection construction (3) and this is an instanton.

\

For another point of view on this we consider  what could be called the
\lq\lq abstract instanton equations'' for a quadruple $S_{0}, S_{1}, S_{2}, S_{3}$ of skew adjoint operators on some vector space. These are the three equations
\begin{equation} [S_{0}, S_{i}] + [S_{j},S_{k}]=0 \end{equation}

where $(ijk)$ runs over cyclic permutations of $(123)$. If we write
$$  D_{1}= S_{0}+ \sqrt{-1} S_{1}\ \ \ , D_{2}= S_{2}+\sqrt{-1} S_{3} , $$
the equations can be written as
\begin{equation} [D_{1}, D_{2}]=0 \ \  ,\ \   [D_{1}, D_{1}^{*}] + [D_{2}, D_{2}^{*}]=0 \end{equation}
If we have a connection on a bundle over $\bR^{4}$ or a flat torus $\bR^{4}/\Lambda$ and let $S_{i}$ be the covariant derivatives $\nabla_{i}$ in the co-ordinate directions  then (10) is just the instanton equation. The operators $D_{1}, D_{2}$ are the components of the coupled $\overline{\partial}$ operator  when we  identify $\bR^{4}$ with $\bC^{2}$ and the commutator equation $[D_{1}, D_{2}]=0$ is the integrability condition for this $\overline{\partial}$ operator to induce a holomorphic structure on the bundle. But in adopting the form (11) we break the natural symmetry of the equations (10): there is a 2-sphere family of equivalent descriptions
corresponding to the $2$-sphere of compatible complex structures. The resulting family of holomorphic bundles is essentially the same as the holomorphic bundle on twistor space defined by Ward.

The Mukai transform on a 4-torus  can be seen as a transform from one solution of the equations (10) to another. The ADHM construction fits into almost the same picture in that if we write $S_{i}=T_{i}$, so that $D_{i}=\alpha_{i}$ the ADHM equations (8)(9) are the same as (11) except for the addition of terms involving $P$ and $Q$. These can be thought of as correction terms due to the non-compactness of $\bR^{4}$. Starting with Nahm \cite{kn:Nahm},  many other examples of such transforms have been found, involving solutions of the instanton equation on $\bR^{4}$ invariant under groups of isometries (see the survey of Jardim \cite{kn:Jardim}).  Nahm studied \lq\lq monopoles'',
which are solutions invariant under  translation in one variable, and defined a transform which related these to solutions of the ODE system
\begin{equation}  \frac{d T_{i}}{ds} = -[ T_{j}, T_{k}] , \end{equation}
(since called Nahm's equations). Here $T_{i}$ for $i=1,2,3$ are functions of a variable $s$ with values in  skew-adjoint matrices. The equations can be written in the form (10) by setting $S_{i}=T_{i}$ for $i=1,2,3$, acting as multiplication operators on vector-valued functions, and
$ S_{0}= \frac{d}{ds}$. One can also think of the data $T_{i}$ as defining a connection over $\bR^{4}$ invariant under translation in three directions and then (12) are just the instanton equations. Nahm's construction follow the same general pattern as those of  ADHM and Mukai,  using solutions of coupled Dirac equations. There are links with both twistor theory and integrable systems which were developed by Hitchin \cite{kn:Hitchin2}.
 The commutator equation $[D_{1}, D_{2}]=0$ is
 $$   \frac{d}{ds} (T_{2}+ \sqrt{-1}T_{3}) = - \sqrt{-1}[ T_{1}, T_{2}+ \sqrt{-1} T_{3}], $$
 which implies that the spectrum of $T_{2}+ \sqrt{-1}T_{3}$ is independent of $s$. We have a $2$-sphere family of similar equations and the resulting family of spectra forms a \lq\lq spectral curve'': a Riemann surface $\Sigma$ with a branched covering over $S^{2}$.  For each $s$ the eigenspaces define a line bundle over $\Sigma$ and Hitchin showed that Nahm's equations correspond
to  linear motion  on the Jacobian of $\Sigma$.


\begin{thebibliography}{99}
\bibitem{kn:AtiyahP}  M. Atiyah {\em  Geometry of Yang-Mills fields} Lezioni Fermiane Acad. Nazionale dei Lincei,  Scuola Normale Sup. Pisa (1979)
\bibitem{kn:Atiyah}  M. Atiyah {\em  Collected works. Vol. 5. Gauge theories} Oxford Science Publications. Oxford UP 1988 
\bibitem{kn:ADHM}  M. Atiyah, V. Drinfeld, N. Hitchin and Y. Manin {\em   Construction of instantons} Physics Letters A 65 (1978) 185-187.
\bibitem{kn:AHS}  M. Atiyah, N. Hitchin and I. Singer {\em   Self-duality in four-dimensional Riemannian geometry} Proc. Roy. Soc. Lond. A 362 425-61 (1978)
\bibitem{kn:BarthHulek} W. Barth and K. Hulek {\em   Monads and moduli of vector bundles} Manuscripta Math. 25  323-347 (1978)
\bibitem{kn:Beilinson}  A. Belinson {\em Coherent sheaves on $\bP^{n}$ and problems in linear algebra} (Russian) Funktsional. Anal. i Prilozhen. 12 (1978) 68-69. 
\bibitem{kn:Buchdahl} N. Buchdahl {\em Instantons on $\bC\bP^{2}$} J. Differential Geom. 24 (1986) 19-52
\bibitem{kn:CG} E. Corrigan and P. Goddard  Corrigan, E.; Goddard, P.{\em  Construction of instanton and monopole solutions and reciprocity} Ann. Physics 154 (1984) 253-279
\bibitem{kn:DK} S. Donaldson and P. Kronheimer {\em The geometry of four-manifolds} Oxford UP 1990
\bibitem{kn:GH} P. Griffiths and J. Harris {\em Principles of algebraic geometry}
Wiley 1978
\bibitem{kn:Hitchin1} N. Hitchin {\em Construction of monopoles} Commun. Math. Physics 89 (1983) 145-190
\bibitem{kn:Hitchin2} N. Hitchin {\em Michael Atiyah: Geometry and Physics} To appear in {\it Literature and history of Mathematical
Science}   International Press
\bibitem{kn:Horrocks}   G. Horrocks{\em  Vector bundles on the punctured spectrum of a local ring} Proc. London Math. Soc.  14  689-713 (1964)
\bibitem{kn:Jardim} M. Jardim {\em A survey on Nahm transform} Jour. Geometry and Physics 52 (2004) 313-327
\bibitem{kn:Mukai} S. Mukai {\em  Duality between $D(X)$ and $D(\hat{X})$ with its application to Picard sheaves} Nagoya Math. J. 81 (1981), 153-175
\bibitem{kn:Nahm} W. Nahm {\em The construction of all self-dual monopoles by the ADHM method} in N. Craigie {\it et al}  Eds. Monopoles in quantum theory World Scientific Singapore 1982
\bibitem{kn:Ward}  R. Ward   {\em On self-dual gauge fields} Phys. Lett. A 61 (1977) 81-82 
\bibitem{kn:Yau} S-T. Yau {\em Shiing-Shen Chern: a great geometer of the twentieth century} To appear in {\it Literature and history of Mathematical Science}   International Press

\end{thebibliography}
\end{document}